\documentclass{amsart}
\usepackage{amssymb}
\usepackage{amsmath}
\usepackage{amsthm} 
\usepackage{hyperref}
\usepackage{xcolor}
\hypersetup{
    colorlinks,
    linkcolor={black!50!black},
    citecolor={blue!50!black},
    urlcolor={blue!80!black}
}
\usepackage[sort&compress,noabbrev, nameinlink]{cleveref}
\usepackage[numbers]{natbib}
\usepackage{url}
\usepackage{enumitem}

\newtheorem{thm}{Theorem}[section]
\newtheorem{lem}[thm]{Lemma}
\newtheorem{cor}[thm]{Corollary}

\newtheorem{defn}{Definition}[section]
\newtheorem{eg}[defn]{Example}
\newtheorem{rmk}[defn]{Remark}
\newtheorem{nott}[defn]{Notation}

\usepackage[normalem]{ulem}
\usepackage{graphicx}

\begin{document}

\title{On m-metric spaces and fixed point theorems.}
\author{Samer Assaf}
\address{Department of Mathematics and Natural Sciences
\\American University of Kuwait
\\Salem Mubarak Street,
Salmiya, Kuwait}
\email{sassaf@auk.edu.kw}
\begin{abstract}
    In this paper we make some observations concerning M-metric spaces and  point out some discrepancies in the proofs found in the literature. To remedy this, we propose a new topological construction and  prove that it is in fact a generalization of a partial metric space. Then, using this construction, we present our main theorem having as its corollaries the fixed point theorems found in previous publications.\end{abstract}
\maketitle

\tableofcontents

\section{Introduction}
\label{S1}
 In 2014  Asadi et al.\cite{Asa2014} proposed the $M-$metric, an intended generalization of a partial metric. In their paper, as we demonstrate in our \Cref{E2b}, the proof of  [\cite{Asa2014}: Lemma 2.5] does not hold. Although [\cite{Asa2014}: Lemma 2.5]  is a small lemma, its assertion was crucial to the proof  of their main theorems; [\cite{Asa2014}: Theorem 3.1] and [\cite{Asa2014}: Theorem 3.2]. Our main concern in their approach lies in the open balls they proposed. We go more in depth on the subject in \Cref{S4}.
\\\indent In this paper, we organize our work in the following manner:
\\ In \Cref{S2}, we introduce the $M-$metric presented in \cite{Asa2014} and generalize it  to allow  negative values. We also present examples that show why some  assumptions proposed in \cite{Asa2014}, including [\cite{Asa2014}: Lemma 2.5], were not accurate. 
\\In \Cref{S3}, we present the partial metric found in \cite{Ass20151,Mat1994,One1996}. We also show how to induce a partial metric from an $M-$metric. The purpous of this section is to put in perspective the generalization from a partial metric to an  $M-$metric.
\\In \Cref{S4}, we discuss why Asadi et. al.'s \cite{Asa2014} proposed open balls were not optimal.  We then present an alternative definition of open balls and discuss the resulting topology. 
\\In \Cref{S5}, we use the topology presented in \Cref{S4} to define  limits and  Cauchy-like sequences in $M-$metric spaces. We then present some of their topological properties. 
\\In \Cref{S6}, we present contractive criteria on functions allowing them to generate Cauchy-like sequences.
\\In \Cref{S7}, we discuss weak orbital continuity, non expansiveness and the lower bound of a space. Those properties are needed for our main theorem.
\\ Finally, in \Cref{S8}, we introduce our main theorem:\\ \textbf{\Cref{T8a}:} \textit{Let $(X,\sigma)$ be an  $M-$metric space with $x_o \in X$. 
Let $f:X\to X$ be a function such that $f$ is $r-$Cauchy at $x_0$ with special limit $a\in X$. Further assume at least one of the following conditions holds:
\begin{enumerate}
\item $f$ is weakly orbitally continuous at $x_0$ and  non-expansive.
\item $f$ is weakly orbitally continuous at $x_0$ and $(X, \sigma)$ is bounded below by $\sigma(f(a), f(a))$.
\item $f$ is non-expansive and $(X, \sigma)$ is bounded below by $\sigma(a, a)$.
\end{enumerate}
Then, $a$ is a fixed point of $f$.}
\\We then use \Cref{T8a} to present a valid proof of [\cite{Asa2014}: Theorem 3.1] and [\cite{Asa2014}: Theorem 3.2]. 
\section{$M-$metric}\label{S2}
\begin{defn}\label{D2a}  Consider a set $X$ and a function $\sigma: X\times X \to\mathbb{R}$. Let
$$m_{x,y}=\min\{\sigma(x,x),\sigma(y,y)\}$$
and
$$M_{x,y}=\min\{\sigma(x,x),\sigma(y,y)\}.$$ 
We say that $\sigma$ is an \textbf{\uline{m-metric}} on $X$ if it satisfies the following axioms:
\\For all $x,y,z\in X,$
\\($\sigma-$lbnd):  $m_{x,y}\le \sigma(x,y)$.
\\($\sigma-$sym):  $\sigma(x, y) = \sigma(y, x)$.
\\($\sigma-$sep):  $\sigma(x,x)=\sigma(x,y)=\sigma(y,y) \iff x=y$.
\\($\sigma-$inq): $\sigma(x, y)-m_{x,y} \le \sigma(x, z)-m_{x,z} + \sigma(z, y)-m_{z,y}.$
\end{defn}
It is important to notice that for all $x,y \in X$, 
$$m_{x,x}=M_{x,x}=\sigma(x,x).$$
\begin{rmk}\label{R2a} In \cite{Asa2014}, the $M-$metric was restricted to having non negative values. In \Cref{D2a}, we remove that restriction to expand on the generalization. The reader may notice that ($\sigma-$lbnd) is redundant since it can be obtained from ($\sigma-$inq). Nevertheless, we chose to state ($\sigma-$lbnd) not only to be consistent with \cite{Asa2014}, but mainly because we use it often enough  to warranty giving it its own name.  
\end{rmk}
 In [\cite{Asa2014}, Example 1.2], the function
\begin{center} $\sigma^\star(x,y)=\sigma(x,y)-m_{x,y}$ for $x\ne y$ and $\sigma^\star(x,x)=0$ \end{center}
was proposed to be a metric.  We present a counterexample in \Cref{E2a}.
\begin{eg}
\label{E2a} Let $\sigma$ be an $M-$metric on the set $X=\{a,b\}$ defined as $$\sigma(a,a)=\sigma(a,b)=\sigma(b,a)=1 \textup{ and } \sigma(b,b)=2.$$ Hence, $m_{a,b}=\min\{\sigma(a,a),\sigma(b,b)\}=1$. Therefore, $$\sigma^\star(a,b)=\sigma(a,b)-m_{a,b}=1-1=0 \textup{ but } a\ne b.$$
Since $\sigma(a,b)-m_{a,b}=0$ with $a\ne b$, $\sigma^\star$  fails to satisfy the metric separation axiom.
\end{eg}
One of the basic ideas behind the $M-$metric is ($\sigma-$lbnd). This axiom ensures that $m_{x,y}=\min\{\sigma(x,x),\sigma(y,y)\}$ is bounded above by $\sigma(x,y)$. Alternatively, $M_{x,y}=\max\{\sigma(x,x),\sigma(y,y)\}$ remains free from any restrictions. This idea is reinforced by ($\sigma-$inq) which cannot be used to bound $M_{x,y}$. That is why, the claim in [\cite{Asa2014}, Lemma 2.5B]:  
$$\lim_{n \to+ \infty}\sigma(x_n,x_{n-1})=0 \Rightarrow\lim_{n \to+ \infty}\sigma(x_n,x_{n})=0,$$
is incorrect. We present the counterexample below.
\begin{eg}\label{E2b} Consider the sequence $\{x_n\}_{n\in\mathbb{N}}$ on a set $X=\{a,b\}$ such that
$$x_n=\begin{cases} a & \textup{if n is odd} \\
b & \textup{ if n is even} \\
\end{cases}$$
Let $\sigma$ be an $M-$metric on $X$ defined by $$\sigma(a,a)=\sigma(a,b)=\sigma(b,a)=0 \textup{ and } \sigma(b,b)=1.$$ 
For example $$\sigma(x_1,x_3)=\sigma(x_1,x_1)=\sigma(a,a)=0=\sigma(a,b)=\sigma(x_1,x_2)=\sigma(x_2,x_3).$$
On the other hand, $$\sigma(x_2,x_4)=\sigma(x_2,x_2)=\sigma(b,b)=1.$$
Therefore, for all $i$,
$$\sigma(x_{i+2},x_{i+1})\le c\sigma(x_{i+1},x_i)$$
satisfying the requirement of [\cite{Asa2014}: Lemma2.5]. It is clear that
$$\lim_{n \to+ \infty}\sigma(x_n,x_{n-1})=0\textup{ while } \lim_{n \to+ \infty}\sigma(x_n,x_n)\textup{ does not exist}$$  as it alternates between $0$ and $1$.
\end{eg}
[\cite{Asa2014}, Lemma 2.5B] is crucial for [\cite{Asa2014}, Lemma 2.5D] and the fixed point theorems presented  in \cite{Asa2014} to hold. Therefore, the techniques used to prove the theorems found in \cite{Asa2014} are no longer valid.
\section{Partial metric}
\label{S3}
As mentioned in \Cref{S1}, the $M-$metric was proposed to generalize the partial metric. In \Cref{D2a}, we expanded on the definition of an $M-$metric space found in \cite{Asa2014} to \ allow negative values. Hence, our $M-$metric is a generalization of the partial metric as defined by O'Neill \cite{One1996}.
\begin{defn}\label{D3a}  A partial metric $p$ on a set $X$ is a function $p: X\times X \to\mathbb{R}$  satisfying the following axioms:
\\For all $x,y,z\in X,$
\\($p-$lbnd):  $p(x,x)\le p(x,y)$. 
\\($p-$sym):  $p(x, y) = p(y, x)$.
\\($p-$sep):  $p(x,x)=p(x,y)=p(y,y) \iff x=y$.
\\($p-$inq): $p(x, y) \le p(x, z)+ p(z, y)-p(z,z)$.
\end{defn}
\begin{rmk}\label{R3a}Notice that ($p-$inq) self regulates when $x=z$ i.e. for any arbitrary function $s: X\times X\to \mathbb{R}$,
$$s(x,y)=s(x,x)+s(x,y)-s(x,x).$$
\end{rmk} For examples on partial metrics we refer the reader to \cite{Ass20151,Mat2009,Mat1994,One1996}.
\\\indent  Asadi et al. \cite{Asa2014} showed that any partial metric is an $M-$metric. Another approach to their proof is by using a well known property which we present in  \Cref{L3a}.
\begin{lem}\label{L3a}  Let $(\Gamma,+,\le)$ be and ordered commutative group. Then, for every $\{a,b,c\}\in\Gamma$,
$$\min\{c,a\}+\min\{c,b\}\le c + \min\{a,b\}$$
and
$$c + \max\{a,b\}\le \max\{c,a\}+\max\{c,b\}.$$
\end{lem}
Hence, for an $M-$metric $\sigma$ on a set $X$, for every $x,y,z \in X$,
$$ -\sigma (z,z)\le -m_{x,z} -m_{z,y}+m_{x,y}$$
and 
$$ -M_{x,z} -M_{z,y}+M_{x,y}\le -\sigma (z,z).$$
In \Cref{E3a}, we slightly adapt [\cite{Asa2014}: Example 1.1] to give an $M-$metric that is not a partial metric.
\begin{eg}\label{E3a} Consider the set $X=\mathbb{R}$. Let $\sigma: X \times X \to \mathbb{R}$ be  defined by setting for all $x,y \in X$,
$$\sigma(x,y)=x+y.$$
Then, $\sigma$ is an $M-$metric on $X$.
\end{eg}
\textbf{Proof:} Except for ($\sigma-$inq), the proof of the other axioms is quite straightforward.  Without loss of generality, assume $x \le y \in X$. Then, 
$$m_{x,y}=2x$$
and
$$\sigma(x,y)-m_{x,y}=x+y-2x=y-x.$$
Let $z \in X$ such that:
\\\uline{Case 1}: $x \le y \le z$. Then, 
$$[\sigma(x, z)-m_{x,z} ]+ [\sigma(z, y)-m_{z,y}]=[z-x]+[z-y]$$
$$\ge z-x\ge y-x=\sigma(x,y)-m_{x,y}.$$
\\\uline{Case 2}: $x \le z \le y$. Then,
$$[\sigma(x, z)-m_{x,z} ]+ [\sigma(z, y)-m_{z,y}]=[z-x]+[y-z]$$
$$= y-x=\sigma(x,y)-m_{x,y}.$$
\\\uline{Case 3}: $z \le x \le y$. Then,
$$[\sigma(x, z)-m_{x,z} ]+ [\sigma(z, y)-m_{z,y}]=[x-z]+[y-z]$$
$$\ge y-z \ge y-x =\sigma(x,y)-m_{x,y}.$$
Clearly, if $x<y\in X$ then
$$2x<x+y<2y$$
i.e. $$m_{x,y}<\sigma(x,y)<\sigma(y,y).$$
Hence, $\sigma$ is not a partial metric.$\hspace{6ex}\square$
\\\indent We now show that an $M-$metric on a set $X$ induces a partial metric on $X$.
\begin{thm}\label{T3a} Let $\sigma$ be an $M-$metric on a set $X$. As in \Cref{D2a}, we denote $M_{x,y}=\max\{\sigma(x,x),\sigma(y,y)\}$ and $m_{x,y}=\min\{\sigma(x,x),\sigma(y,y)\}.$  For  $x,y \in X$, let
 $$p^{\sigma}(x,y)=\sigma(x,y)+M_{x,y}-m_{x,y}.$$
 Then, $p^{\sigma}$ is a partial metric on $X$.
\end{thm}
\textbf{Proof:} For all $x\in X$, $p^{\sigma}(x,x)=\sigma(x,x)+M_{x,x}-m_{x,x}=\sigma(x,x)$. The proof of  ($p-$sym) is trivial.
\\($p-$lbnd): For all $x,y\in X$, from ($\sigma-$lbnd) we have $\sigma(x,y)-m_{x,y}\ge 0$. Hence, 
$$p^{\sigma}(x,x)=\sigma(x,x)\le M_{x,y}\le M_{x,y}+\sigma(x,y)-m_{x,y}=p^{\sigma}(x,y).$$
($p-$sep): Assume that $p^{\sigma}(x,x)=p^{\sigma}(x,y)=p^{\sigma}(y,y)$. Then, 
$$\sigma(x,x)=p^{\sigma}(x,x)=p^{\sigma}(y,y)=\sigma(y,y)$$
i.e.
$$\sigma(x,x)=\sigma(y,y)=m_{x,y}=M_{x,y}.$$
 Therefore, 
 $$p^{\sigma}(x,y)=\sigma(x,y)+M_{x,y}-m_{x,y}=\sigma(x,y).$$ Hence, $\sigma(x,x)=\sigma(x,y)=\sigma(y,y)$ and, therefore, by ($\sigma-$sep) $x=y$.
\\($p-$inq): For all $x,y,z \in X$,
$$p^{\sigma}(x,y)=[\sigma(x,y)-m_{x,y}]+M_{x,y}$$
by ($\sigma-$inq)
$$\le [\sigma(x,z)-m_{x,z}+\sigma(z,y)-m_{z,y}]+M_{x,y}$$
$$=p^{\sigma}(x,z)+p^{\sigma}(z,y)-M_{x,z}-M_{z,y}+M_{x,y}$$
by \Cref{L3a}
$$\le p^{\sigma}(x,z)+p^{\sigma}(z,y)-\sigma(z,z)$$
$$=p^{\sigma}(x,z)+p^{\sigma}(z,y)-p^{\sigma}(z,z).\hspace{6ex}\square$$
\indent As shown in \Cref{E2a},  given $\sigma$ an $M-$metric on a set $X$,
$$\sigma^\star(x,y)=\sigma(x,y)-m_{x,y}$$ need not  be a metric. However, we can guarantee that  $\sigma^\star$ is a metric in the special case where $\sigma$ is a partial metric. 
\begin{lem}\label{L3b} Let $\sigma$ be a partial metric on a set $X$. For all $x,y \in X$, let
$$\sigma^\star(x,y)=\sigma(x,y)-m_{x,y}.$$
Then, $\sigma^\star$ is a metric on $X$.
\end{lem}
\textbf{Proof:} The major issue in \Cref{E2a} was the metric separation axiom. Since $\sigma$ is a partial metric, and by ($p-$lbnd), for all $x,y \in X$
$$ \sigma(x,y)-\sigma(x,x)\ge 0.$$
Therefore,
$$0\le \sigma(x,y)-\sigma(x,x)\le \sigma(x,y)-m_{x,y}=\sigma^\star(x,y).$$
Hence, if $\sigma^\star(x,y) = 0$ then $\sigma(x,x)=\sigma(x,y)=\sigma(y,y)$ and by ($p-$sep) $x=y$. The rest of the axioms are straightforward and easy to check.  
\section{Topology}
\label{S4}
Let $\sigma$ be an  $M-$metric on a set $X$. For every $x \in X$ and $\epsilon>0$,  Asadi et al. \cite{Asa2014}
defined an $A-$open ball as
$$B^A_\epsilon(x)=\{y\in X \vert \sigma^\star(x,y)=\sigma(x,y)-m_{x,y}<\epsilon\}.$$
Therefore, for the special case of $\sigma$ being a partial metric and from \Cref{L3b}, the $A-$open balls span a metric space. Moreover, in [\cite{Asa2014}: Theorem 2.1], the authors stated that the topology generated by the $A-$open balls is not Hausdorff. This is a faux pas since a metric  is an $M-$metric. 
\\\indent  On the other hand, given a partial metric $p$ on a set $X$, Matthews \cite{Mat1994}  defined the $p-$open ball as$$B^p_\epsilon(x)=\{y \in X \vert p(x,y)-p(x,x)<\epsilon\}.$$
Matthews \cite{Mat1994}  also showed that the $p-$open balls span a $T_0$ topology that need not be $T_1$. We will call the $p-$open balls  the standard partial metric balls. In \cite{Ass20151}, we showed that the standard partial metric balls still work when allowing the partial metric to have negative values i.e. taken in the sense of O'Neill \cite{One1996}.
\begin{lem}\label{L33b} Let $\sigma$ be a partial metric on a set $X$. $\mathcal{T}_A$, the topology generated by $B^A_\epsilon$ balls, is finer than $\mathcal{T}_{\sigma^s}$, the standard partial metric topology generated by the $B^p_\epsilon$ balls.
\end{lem}
\textbf{Proof:} For all $x\in X$, for each $y\in B^p_\epsilon(x)$ i.e. $\sigma(x,y)-\sigma(x,x)<\epsilon$, let $$\delta=\epsilon-\sigma(x,y)+\sigma(x,x)>0.$$ We show that $B^A_\delta(y)\subseteq B^p_\epsilon(x)$. If $z \in B^A_\delta (y)$ i.e. $\sigma(y,z)-m_{y,z}<\delta$, and using ($p-$inq) and $-\sigma(y,y)\le -m_{y,z}$, we get  
$$\sigma(x,z)-\sigma(x,x)$$
$$\le \sigma(x,y)+\sigma(y,z)-\sigma(y,y)-\sigma(x,x)$$
$$\le \sigma(x,y)-\sigma(x,x)+\sigma(y,z)-m_{y,z}$$
$$<\sigma(x,y)-\sigma(x,x)+\delta=\epsilon.$$ And, hence, $B^A_\delta(y)\subseteq B^p_\epsilon(x)$ i.e. $\mathcal{T}_A$ is finer than $\mathcal{T}_{\sigma^s}$.
\\
\\\indent \Cref{L33b} shows why $B^A_\epsilon(x)$ is not an optimal generalization of $B^{p}_\epsilon (x)$. In  the special case where $\sigma$ is a partial metric, $B^A_\epsilon(x)$ becomes a metric ball. Hence, in that case, the topology generated by $B^A_\epsilon(x)$ is much finer than the one generated by $B^p_\epsilon (x)$.
\\\indent  If an $M-$metric theory is to be developed as a generalization of the partial metric one, the topology proposed should  not be finer  than the standard partial metric topology.
The $M-$open balls presented below accomplish just that. We define them and show that the collection of  $M-$open balls form a basis.  
\begin{defn}
\label{D4a} Let $\sigma$ be an $M-$metric on a set $X$. For every $x \in  X$ and $\epsilon>0$, the \uline{\textbf{$M-$open ball}} around $x$ of radius $\epsilon$ is
$$B^\sigma_\epsilon(x)=\{y\in X \vert \sigma(x,y) +\sigma(y,y)-m_{x,y}-\sigma(x,x)<\epsilon\}.$$
\end{defn}
\begin{rmk}\label{R4a} We notice that $\sigma(x,x) -\sigma(x,x)+\sigma(x,x)-m_{x,x}=0$ i.e. for every $\epsilon>0$, $x \in B^\sigma_\epsilon(x)$. Additionally, if  for some $x,y \in X$
$$m_{x,y}=\sigma(y,y)\le\sigma(x,y)\le\sigma(x,x),$$ 
then
$$ \sigma(x,y) +\sigma(y,y)-m_{x,y}-\sigma(x,x)$$
 $$=[\sigma(x,y) -\sigma(x,x)]+[\sigma(y,y)-m_{x,y}]$$
 $$=\sigma(x,y) -\sigma(x,x)\le0.$$
Hence, every $M-$open ball centered at $x$ contains $y$.
\end{rmk}
\begin{lem} \label{L4a} Let $\sigma$ be an $M-$metric on a set $X$. The collection of all $M-$open balls on $X$, $\mathcal{B}^\sigma =\{B^\sigma_\epsilon (x)\}_{x\in X}^{\epsilon >0}$ forms a basis on $X$.
\end{lem}
\textbf{Proof:} For every $x\in X$ and $\epsilon>0$, let $y \in B^\sigma_\epsilon(x)$. Then,$$ \sigma(x,y) +\sigma(y,y)-m_{x,y}-\sigma(x,x)<\epsilon.$$
Take $$\delta = \epsilon - \sigma(x,y) -\sigma(y,y)+m_{x,y}+\sigma(x,x)>0.\hspace{4ex}(\star)$$
We claim that $B^\sigma_\delta(y)\subseteq B^\sigma_\epsilon(x).$ If $z \in B^\sigma_\delta(y)$, then 
$$\sigma(y,z) +\sigma(z,z)-m_{y,z}-\sigma(y,y)<\delta.\hspace{3ex}(\oplus)$$
Hence, by ($M-$inq) (see  \Cref{D2a})
$$\sigma(x,z)+\sigma(z,z)-m_{x,z} -\sigma(x,x)=[\sigma(x,z)-m_{x,z}]+\sigma(z,z)-\sigma(x,x) $$
$$\le [\sigma(x,y)-m_{x,y}+\sigma(y,z)-m_{y,z} ]+\sigma(z,z)-\sigma(x,x)$$
by adding and subtracting $\sigma(y,y)$ we get
$$=[\sigma(x,y) +\sigma(y,y)-m_{x,y}-\sigma(x,x)]+[\sigma(y,z) +\sigma(z,z)-m_{z,z}-\sigma(y,y)].$$
By $(\oplus)$ and $(\star)$, we get
$$\sigma(x,z)+\sigma(z,z)-m_{x,z} -\sigma(x,x)<\sigma(x,y) +\sigma(y,y)-m_{x,y}-\sigma(x,x)+\delta=\epsilon.$$
Therefore,  $B^\sigma_\delta(y)\subseteq B^\sigma_\epsilon(x)$ and $\mathcal{B}^\sigma$ is a basis on $X$.$\hspace{6ex}\square$
\begin{nott}\label{N4a} Given an $M-$metric $\sigma$ on a set $X$. We denote by:
\\  \uline{$\mathcal{T}_\sigma$:} The topology generated by the $M-$open balls $$B^\sigma_\epsilon(x)=\{y\in X \vert \sigma(x,y) +\sigma(y,y)-m_{x,y}-\sigma(x,x)<\epsilon\}.$$ 
\uline{$\mathcal{T}_{p^\sigma}$:} The standard partial metric topology spanned by the $p-$open balls
 $$B^{p^\sigma}_\epsilon(x)=\{y \in X \vert p^\sigma(x,y)-p^\sigma(x,x)<\epsilon\},$$ where $p^\sigma$ is the induced partial metric defined in \Cref{T3a}. I.e. 
$$B^{p^\sigma}_\epsilon(x)=\{y \in X \vert \sigma(x,y)+M_{x,y}-m_{x,y}-\sigma(x,x)<\epsilon\}.$$
In the special case where $\sigma$ is a partial metric, we denote by
\\ \uline{$\mathcal{T}_{\sigma_s}$:} The standard partial metric  topology generated by the $p-$open balls
 $$B^{\sigma_s}_\epsilon(x)=\{y \in X \vert \sigma(x,y)-\sigma(x,x)<\epsilon\}.$$
 \end{nott}
We now move to comparing the topologies defined in \Cref{N4a}.  Given an $M-$metric $\sigma$ on a set $X$, we show in \Cref{L4b}  that $\mathcal{T}_\sigma$ is coarser than $\mathcal{T}_{p^\sigma}$. In the special case where $\sigma$ is a partial metric, we show in \Cref{L4c} that $\mathcal{T}_{\sigma_s}=\mathcal{T}_\sigma$.
\Cref{L4b} and \Cref{L4c} sheds light as to why we consider $\mathcal{T}_\sigma$ to be a proper generalization of $\mathcal{T}_{\sigma_s}$.
\begin{lem}\label{L4b}
 Let $\sigma$ be an $M-$metric on a set $X$. Then, $\mathcal{T}_\sigma$ is coarser than $\mathcal{T}_{p^\sigma}$. 
\end{lem}
\textbf{Proof:} Consider the $M-$open ball $B^\sigma_\epsilon (x)$.  Let $y\in B^{p^\sigma}_\epsilon(x)$. Then, from \Cref{N4a},
$$ \sigma(x,y)+M_{x,y}-m_{x,y}-\sigma(x,x)<\epsilon.$$
Hence, $$\sigma(x,y)+\sigma(y,y)-m_{x,y}-\sigma(x,x)\le\sigma(x,y)+M_{x,y}-m_{x,y}-\sigma(x,x)<\epsilon.$$
Therefore, $B^{p^\sigma}_\epsilon(x)\subseteq B^\sigma_\epsilon(x)$ and, hence, $\mathcal{T}_\sigma \subseteq \mathcal{T}_{p^\sigma}.\hspace{6ex}\square$
\begin{lem}\label{L4c} Let $\sigma$ be a partial metric on a set $X$. Then, $\mathcal{T}_{\sigma_s}=\mathcal{T}_\sigma$.
\end{lem}
\textbf{Proof:} Consider the standard  $p-$open ball $B^{\sigma_s}_\epsilon (x)$. Let $y\in B^{\sigma}_\epsilon(x)$. Then, $ \sigma(x,y) +\sigma(y,y)-m_{x,y}-\sigma(x,x)<\epsilon$.
$$\sigma(x,y)-\sigma(x,x)\le \sigma(x,y)-\sigma(x,x)+\sigma(y,y)-m_{x,y}<\epsilon.$$
Hence, $B^{\sigma}_\epsilon(x)\subseteq B^{\sigma_s}_\epsilon(x)$. Therefore, $\mathcal{T}_{\sigma_s} \subseteq \mathcal{T}_{\sigma}$.
\vspace{2ex}
\\\indent Conversely, let $y \in B^{\sigma_s}_{\frac{\epsilon}{2}}(x)$. Then, $\sigma(x,y)-\sigma(x,x)<\frac{\epsilon}{2}.$
\\Case 1: If $m_{x,y}=\sigma(x,x)\le\sigma(y,y)$ then from ($p-$lbnd) in \Cref{D3a},
 $$\sigma(y,y)\le \sigma(x,y) \textup{ and } m_{x,y}=\sigma(x,x).$$ 
 Hence,
 $$\sigma(x,y)+\sigma(y,y)-m_{x,y}-\sigma(x,x)\le \sigma(x,y)+\sigma(x,y)-\sigma(x,x)-\sigma(x,x)<2(\frac{\epsilon}{2})=\epsilon.$$
Case 2: $m_{x,y}=\sigma(y,y)\le\sigma(x,x)$ then, $\sigma(y,y)-m_{x,y}=0$ and, hence,
$$ \sigma(x,y)+\sigma(y,y)-m_{x,y}-\sigma(x,x)=\sigma(x,y)-\sigma(x,x)<\frac{\epsilon}{2}<\epsilon.$$
Therefore, $B^{\sigma_s}_{\frac{\epsilon}{2}}(x)\subseteq B^\sigma_\epsilon(x)$ and, hence, $\mathcal{T}_\sigma=\mathcal{T}_{\sigma_s}$. $\hspace{6ex}\square$ 
\begin{nott}\label{N4b} Let $\sigma$ be an $M-$metric on a set $X$. We denote by:
\\The $M-$metric space $(X,\sigma)=(X,\mathcal{T}_\sigma)$.
\\The partial metric space $(X, p^\sigma)=(X,\mathcal{T}_{p^\sigma})$.
\\We remind the reader that if $\sigma$ is a partial metric, then the standard partial metric space $(X,\mathcal{T}_{\sigma_s})=(X,\mathcal{T}_\sigma)=(X,\sigma)$.
\end{nott}
All our work would be useless if for every $M-$metric $\sigma$, $\mathcal{T}_\sigma= \mathcal{T}_{p^\sigma}.$ We use the $M-$metric defined in \Cref{E3a} to give an example where $\mathcal{T}_\sigma \subsetneq \mathcal{T}_{p^\sigma}$.
\begin{eg}\label{E4a} Let $\sigma$ be an $M-$metric on $X=\mathbb{R}$ as defined in \Cref{E3a} by
$$\sigma(x,y)=x+y.$$
Then, $\mathcal{T}_\sigma \subsetneq \mathcal{T}_{p^\sigma}$.
\end{eg}
\textbf{Proof:} We remind our reader  that if $x \le y \in X$ then
$$\sigma(y,y)=2y \textup{ and }\sigma(x,x)=m_{x,y}=2x.$$
From \Cref{N4a}, $\mathcal{T}_\sigma$ is generated by the $M-$open balls
$$B^\sigma_\epsilon(x)=\{y\in X \vert \sigma(x,y) +\sigma(y,y)-m_{x,y}-\sigma(x,x)<\epsilon\}.$$
If $y \le x$ then $\sigma(x,y) +\sigma(y,y)-m_{x,y}-\sigma(x,x)=y-x\le 0<\epsilon$.
\\If $x<y$ then $\sigma(x,y) +\sigma(y,y)-m_{x,y}-\sigma(x,x)=3(y-x)$. Therefore, $$ B^\sigma_\epsilon(x)=(-\infty,x+\frac{\epsilon}{3}).$$
Again from \Cref{N4a}, $\mathcal{T}_{p^\sigma}$ is generated by the $p-$ open balls
$$B^{p^\sigma}_\epsilon(x)=\{y \in X \vert \sigma(x,y)+M_{x,y}-m_{x,y}-\sigma(x,x)<\epsilon\}.$$
If $y \le x$ then  $\sigma(x,y)+M_{x,y}-m_{x,y}-\sigma(x,x)=x-y$.
\\If $x<y$ then $\sigma(x,y)+M_{x,y}-m_{x,y}-\sigma(x,x)=3(y-x)$. Therefore,
$$B^{p^\sigma}_\epsilon(x)=(x-\epsilon,x+\frac{\epsilon}{3}).$$
Clearly, $\mathcal{T}_\sigma \subsetneq \mathcal{T}_{p^\sigma}$. $\hspace{6ex}\square$
\begin{lem}
\label{L4d} An $M-$metric space is $T_o$.
\end{lem}
\textbf{Proof:} Let $(X,\sigma)$ be an $M-$metric space with two distinct elements $x, y \in X$. Without loss of generality, we can consider two cases:
\\Case 1: If $\sigma(x,x)=\sigma(y,y)$ then by ($\sigma-$lbnd) and ($\sigma-$sep) (see \Cref{D3a}),  and since $x\ne y$, we have 
$$m_{x,y}=\sigma(x,x)=\sigma(y,y)<\sigma(x,y).$$ 
Hence,
$$\sigma(x,y)+\sigma(y,y)-m_{x,y}-\sigma(x,x)=\sigma(x,y)-\sigma(x,x)>0.$$
Therefore, if $\epsilon=\sigma(x,y)-\sigma(x,x)$ then $y\notin B^\sigma_\epsilon(x).$
\\Case 2: If $\sigma(x,x)<\sigma(y,y)$ then by ($\sigma-$lbnd) 
$$\sigma(x,y)-m_{x,y}\ge 0.$$
Hence, 
$$\sigma(x,y)+\sigma(y,y)-m_{x,y}-\sigma(x,x)\ge \sigma(y,y)-\sigma(x,x)>0.$$
Therefore, if $\epsilon=\sigma(y,y)-\sigma(x,x)$ then $y \notin B^\sigma_\epsilon(x)$.
\\ Since a partial metric is an $M-$metric, we refer the reader to \cite{Mat1994,Ass20151} for examples of $M-$metric spaces that need not be $T_1$.
\section{$r-$Cauchy sequences and Limits}
\label{S5}
We begin this section by defining a Cauchy-like sequence. We use the same approach found in \cite{Ass20151,Mat2009} and apply it to the $M-$metric case.
\begin{defn}\label{D5a}
 Let $(X,\sigma)$ be an $M-$metric space and $r$ a real number. A sequence $\{x_i\}_{i\in\mathbb{N}}$ in $X$ is said to be \textbf{\uline{$r-$Cauchy}} iff
$$\lim_{i,j\to +\infty}\sigma(x_i,x_j)=r.$$
$r$ is called the \textbf{\uline{central distance}} of  $\{x_i\}_{i\in\mathbb{N}}$.
\end{defn}
\begin{rmk}
\label{R5a} Alternatively, we could have defined an $r-$Cauchy sequence as:
$$\lim_{i\ne j\to + \infty}\sigma(x_i,x_j)=\lim_{i, j\to +\infty}m_{x_i,x_j}=r.$$ This definition is closer to the one presented in \cite{Asa2014} but since it  has a subsequence $\{y_i\}_{i\in\mathbb{N}}$ such that
$$\lim_{i,j\to \infty}\sigma(y_i,y_j)=r,$$
we find that there is very little point in using a more general definition.
\end{rmk}
\begin{lem}
\label{L5a} Let $\{x_i\}_{i\in\mathbb{N}}$ be an $r-$Cauchy sequence in an $M-$metric space $(X,\sigma)$. Then,
$$\textup{ (a) }\lim_{i\to+ \infty}\sigma(x_i,x_i)=r.$$
$$\textup{ (b) }\lim_{i,j\to+ \infty}m_{x_i,x_j}=r.$$
and
$$\textup{ (c) }\lim_{i,j\to+ \infty}M_{x_i,x_j}=r.$$
\end{lem}
\textbf{Proof:} The proof of \Cref{L5a} is quite straightforward. (a) follows trivially  from \Cref{D5a}. (b) and (c) follow trivially from (a).$\hspace{6ex}\square$
\\\indent The limit is a topological definition which we translate into the  language of $M-$metric spaces.
\begin{defn} 
\label{D5b}Let $\{x_i\}_{i\in\mathbb{N}}$ be an $r-$Cauchy\ sequence in an $M-$metric space $(X,\sigma)$. We say that $a\in X$ is a \textbf{\uline{limit}} of $\{x_i\}_{i\in\mathbb{N}}$ iff
$$\lim_{i\to+ \infty}\sigma(a,x_i)+\sigma(x_i,x_i)-m_{a,x_i}=\sigma(a,a).$$
\end{defn}
The natural question to ask here is: How does $r$ relate to $\sigma(a,a)$? This question has been answered in the partial metric case in \cite{Ass20151}. The result remains the same in an $M-$metric case.
\begin{lem}\label{L5b} Let $\{x_i\}_{i\in\mathbb{N}}$ be an $r-$Cauchy sequence in an $M-$metric space $(X,\sigma)$. If $a$ is a limit of $\{x_i\}_{i\in\mathbb{N}}$ then
$$r\le \sigma(a,a).$$ 
\end{lem}
\textbf{Proof:} By ($\sigma-$lbnd) (see \Cref{D2a}) we know that for each $i$,
$$0\le \sigma(a,x_i)-m_{a,x_i}.$$
Adding $\sigma(x_i,x_i)$ on both sides we get 
$$\sigma(x_i,x_i)\le \sigma(a,x_i)+\sigma(x_i,x_i)-m_{a,x_i}.$$
Now taking the limit of both sides, by \Cref{L5a}(a) and \Cref{D5b}, we get
$$r\le \sigma(a,a).\hspace{6ex}\square$$
\indent The limit of an $r-$Cauchy sequence need not be unique. An example is given for the partial metric case in \cite{Ass20151}.
\\\indent Reading through the partial metric literature, \cite{Ass20151,Kar20112,Kar20131,Mat2009,One1996} to name but a few,  it becomes obvious that a stronger version of a limit is needed. The $M-$metric space, being a generalization of a partial metric space, is no exception.
\begin{defn}
\label{D5c} Let $\{x_i\}_{i\in\mathbb{N}}$ be an $r-$Cauchy sequence in an $M-$metric space $(X,\sigma)$. An element $a\in X$ is a \textbf{\uline{special limit}} of $\{x_i\}_{i\in\mathbb{N}}$ iff 
\begin{center}$a$ is a limit of  $\{x_i\}_{i\in\mathbb{N}}$ and $\sigma(a,a)=r$.
\end{center}
\end{defn}
Unlike a regular limit, a special limit is  unique.
\begin{lem}\label{L5c} Let $\{x_i\}_{i\in\mathbb{N}}$ be an $r-$Cauchy sequence in an $M-$metric space $(X,\sigma)$. If $a$ is a special limit of $\{x_i\}_{i\in\mathbb{N}}$ then $a$ is unique.
\end{lem}
\textbf{Proof:} Let $a$ and $b$ be two  special limits of $\{x_i\}_{i\in\mathbb{N}}$. From \Cref{D5c}, $$r=\sigma(a,a)=\sigma(b,b)=m_{a,b}.$$
From ($\sigma-$lbnd) (see \Cref{D2a}) we know that 
$$r=\sigma(a,a)=m_{a,b}\le \sigma(a,b).$$
Hence,
$$\sigma(a,b)=[\sigma(a,b)-m_{a,b}]+r$$
using ($\sigma-$inq), for all $i$,
$$\le [\sigma(a,x_i)-m_{a,x_i}+\sigma(b,x_i)-m_{b,x_i}]+r.$$
Therefore, by adding and subtracting $\sigma(x_i,x_i)$ we get
$$\sigma(a,b)\le[ \sigma(a,x_i)+\sigma(x_i,x_i)-m_{a,x_i}]+[\sigma(b,x_i)+\sigma(x_i,x_i)-m_{b,x_i}]-2\sigma(x_i,x_i)+r.$$
Since a special limit is also a limit (see \Cref{D5c}) and by \Cref{L5a}:
$$\sigma(a,x_i)+\sigma(x_i,x_i)-m_{a,x_i}\to \sigma(a,a)=r$$,
$$\sigma(b,x_i)+\sigma(x_i,x_i)-m_{b,x_i} \to \sigma(b,b)=r$$
and
$$\sigma(x_i,x_i)\to \sigma(a,a)=r.$$
Hence, 
$$\sigma(a,b)\le r+r-2r+r=r=\sigma(a,a)$$
therefore,
$$\sigma(a,a)=\sigma(a,b)=\sigma(b,b).$$
By ($\sigma-$inq) (see \Cref{D2a}) $a=b.\hspace{6ex}\square$
\begin{lem}
\label{L5d} Let $\{x_i\}_{i\in\mathbb{N}}$ be an $r-$Cauchy sequence in an $M-$metric space $(X,\sigma)$. Let $a$ be the special limit of  $\{x_i\}_{i\in\mathbb{N}}$. Then,
$$\textup{ (a) }\lim_{i\to+ \infty}M_{a,x_i}=\sigma(a,a).$$
$$\textup{ (b) }\lim_{i\to+ \infty}m_{a,x_i}=\sigma(a,a).$$
and
$$\textup{ (c) }\lim_{i\to+ \infty}\sigma(a,x_i)=\sigma(a,a).$$
\end{lem}
\textbf{Proof:} Parts (a) and (b) are straightforward. As for (c), using \Cref{D5b} we get 
$$\lim_{i\to+ \infty}\sigma(a,x_i)+\sigma(x_i,x_i)-m_{a,x_i}=\sigma(a,a).$$
Hence, using \Cref{D5c} and \Cref{L5a}(b) we get
$$\lim_{i\to+ \infty}\sigma(a,x_i)=\lim_{i\to+ \infty}([\sigma(a,x_i)+\sigma(x_i,x_i)-m_{a,x_i}]-\sigma(x_i,x_i)+m_{a,x_i})$$
$$=\sigma(a,a)-\sigma(a,a)+\sigma(a,a)=\sigma(a,a).\hspace{6ex}\square$$
\begin{lem}\label{L5e} Let $\{x_i\}_{i\in\mathbb{N}}$ be an $r-$Cauchy sequence in an $M-$metric space $(X,\sigma)$. If $a$ is a special limit of $\{x_i\}_{i\in\mathbb{N}}$ then for every $y\in X$,
$$\lim_{i\to +\infty} \sigma(y,x_i)-m_{y,x_i}=\sigma (y,a)-m_{y,a}.$$
\end{lem}
\textbf{Proof:} The proof follows directly from ($\sigma-$lbnd) and \Cref{L5d}.
\\\indent What is left is to guarantee the exitance of a special limit. Therefore, we present the notion of completeness in an $M-$metric space.
\begin{defn}\label{D5d} An $M-$metric space $(X,\sigma)$ is said to be \textbf{\uline{complete}} iff for every real number $r$, every $r-$Cauchy sequence in $X$ has a special limit in $X$.
\end{defn}
\section{$r-$Cauchy Functions}\label{S6}
One of the cornerstones of Banach-like fixed point theorems is that the function  $f$ in question has a Cauchy-like  orbit $\{f^i(x_o)\}_{i \in \mathbb{N}}$ for some $x_o\in X$.
  \begin{defn}\label{D7a} Let  $(X,\sigma)$ be an $M-$etric space with $x_o \in X$. Suppose $f:X\to X$ is a function on $X$. We say that $f$ is a \uline{\textbf{$r-$Cauchy function} at $x_o$} if and only if $ \{f^i(x_o)\}_{i \in \mathbb{N}}$ is an $r-$Cauchy sequence in $X$.
\end{defn}
Throughout the literature, different criteria  on a function $f$ were investigated for $f$ to be an $r-$Cauchy function. Many of those cases boil down to two main ones which we  present in \Cref{D7a} and \Cref{D7b}. 
\begin{defn}\label{D7b} Let  $(X,\sigma)$ be an $M-$metric space with $x_o \in X$. Let $f:X\to X$ be a function on $X$. Let $r$ and $0<c<1$ be two real numbers. We say that $f$ is an \uline{\textbf{orbital $c_r-$contraction at $x_o$}} (or $f$ is \uline{\textbf{orbitally $c_r-$contractive at $x_o$}}) if and only if  for all natural numbers $i$,
$$r \le \sigma(f^{i+1}(x_o),f^{i+1}(x_0))\le r+ c^{i}\vert\sigma(f(x_o),x_o)\vert$$
and
$$\sigma(f^{i+2}(x_o),f^{i+1}(x_0) )\le r +c^{i+1}\vert\sigma(f(x_o),x_o)\vert.$$
\end{defn}
\begin{lem}\label{L7a} Let  $(X,\sigma)$ be an $M-$metric space with $x_o \in X$. Let $f:X\to X$ be a function on $X$. Let $r$ and $0<c<1$ be two real numbers. If $f$ is an orbital $c_r-$contraction at $x_o$ then $f$ is an $r-$Cauchy function at $x_o$.
\end{lem}
\textbf{Proof:} The proof is quite straightforward by first showing 
$$\lim_{i \to +\infty}\sigma(x_i,x_i)=r \textup{ and, hence, } \lim_{i,j \to +\infty}m_{x_i,x_j}=r.$$ For a detailed similar proof, we refer the reader to [\cite{Ass20151}: Lemma 6.2].$\hspace{6ex} \square$
\begin{defn}\label{D7c} Let  $(X,\sigma)$ be an $M-$metric space with $x_o \in X$. Let $f:X\to X$ be a function on $X$. Let $r$ be a real number and  $\varphi:[r,+\infty )\subset\mathbb{R}\to [0, + \infty)$  be a non-decreasing function such that
$$\varphi(t)=0 \textup{ if and only if }t=r.$$
We say that $f$ is an \uline{\textbf{orbital $\varphi_{r}-$Contraction at $x_o$}}  if and only if for all $i$ and $j$,
$$r \le \sigma(f^{i+1}(x_o),f^{j+1}(x_o))\le \sigma(f^i(x_o),f^j(x_o))-\varphi (\sigma(f^i(x_o),f^j(x_o)).$$
\end{defn}
The proof of \Cref{L7b} is quite delicate. We will give it in its most explicit form while repeatedly clarifying any ambiguous notation.\begin{lem}\label{L7b} Let  $(X,\sigma)$ be an $M-$metric space with $x_o$ in $X$. Let $f:X\to X$ be a function on $X$. If $f$ is an orbital $\varphi_r$-contraction at $x_o$ then $f$ is an $r-$Cauchy function at $x_o$.
\end{lem}
\textbf{Proof:} Let $x_o \in X$ and suppose $f: X \to X$ is an orbital $\varphi_r-$contraction at $x_o$. Denote  $x_i =f^i(x_o)$.
 To remedy any possible ambiguity, we will be adding parenthesis to differentiate between  $x_{(n_k+1)}$ and $x_{n_{(k+1)}}$ when the need arises. 
 \\\indent \uline{Step 1:} Let $t_i=\sigma(x_{i+1},x_i)$. In this step, we will show that in the topological space $\mathbb{R}$ (endowed with the standard topology) $\{t_i\}_{i\in\mathbb{N}}$ is a Cauchy sequence that converges to $r$.
\\From \Cref{D7c}, for all $i$,
$$\sigma(x_{i+1},x_{i+1})\le \sigma(x_{i},x_{i})-\varphi(\sigma(x_{i},x_{i}))$$
and, hence, $\{\sigma(x_i,x_i)\}_{i\in \mathbb{N}}$ forms a decreasing chain since, for all $i$, 
$\varphi(t_i)\ge 0$ i.e.
$$m_{x_i,x_{i+1}}=\sigma(x_{i+1},x_{i+1}).$$
Moreover,  from($\sigma-$lbnd)  
$$r\le \sigma (x_{i+2},x_{i+2})=m_{x_{i+2},x_{i+1}}\le \sigma(x_{i+2},x_{i+1})=t_{i+1}$$
and
$$t_{i+1}=\sigma(x_{i+2},x_{i+1})\le \sigma(x_{i+1},x_i)-\varphi (\sigma(x_{i+1},x_i))$$
$$=t_i-\varphi(t_i)\le t_i.$$
Hence, for all $i$,
$$r\le t_{i+1}\le t_i$$
i.e. $\{t_i\}_{i\in \mathbb{N}}$ is a non-increasing sequence in $\mathbb{R}$ bounded below by $r$ and, therefore, $\{t_i\}_{i\in \mathbb{N}}$ is a Cauchy sequence in $\mathbb{R}$. Since $\mathbb{R}$ with the standard topology is a complete metric space, $\{t_i\}_{i\in\mathbb{N}}$ has a limit $L$ such that for all $i$,
$$t_i\ge L \ge r$$
and, since $\varphi$ is a non-decreasing function,
$$\varphi(t_i) \ge \varphi(L) \ge \varphi(r)=0$$
i.e.
$$-\varphi(t_i) \le -\varphi(L)\le 0.$$
Hence, by \Cref{D7c}
$$r\le t_{i+1}\le t_i-\varphi(t_i)\le t_i-\varphi(L)$$
$$\le t_{i-1}-\varphi(t_{i-1})-\varphi(L)\le t_{i-1}-2\varphi(L)$$
by induction
$$t_{i+1}\le t_1-i\varphi(L).$$
Assume that $L>r$. Then, by \Cref{D7c}, $\varphi(L)>0$. By taking $i>\frac{t_1-r}{\varphi(L)}$ we get 
$$t_{i+1} \le t_1-i\varphi(L)<t_1-\frac{t_1-r}{\varphi (L)}\varphi(L) =r$$ a contradiction since $t_i \ge r$. Therefore, 
$$\lim_{i \to +\infty }t_i=\lim_{i \to +\infty }\sigma(x_i,x_{i+1})=\lim_{i,j \to + \infty}m_{x_i,x_j}=r.\hspace{6ex}(\bar\bigcirc)$$
\\\indent \uline{Step 2:} We now show that $\{x_i\}_{i\in \mathbb{N}}$ is an $r-$Cauchy sequence by supposing that it is not (a contrapositive approach). 
\\Suppose that  $\{x_i\}_{i\in \mathbb{N}}$ is not an $r-$Cauchy sequence. Since $r\le \sigma(x_i,x_j)$, there exists a positive real number $\delta$ such that for every natural number $N$, there exists  $i,j>N$ where 
$$\sigma(x_i,x_j)\ge r+\delta >r.$$
From step 1, by choosing $N$  big enough, for all  $i>N$,
$$r\le \sigma(x_{i},x_{i})=m_{x_{i-1},x_i}\le \sigma(x_{i},x_{i-1})<r+\delta.$$
Then there exist  $j_1>h_1>N$ such that  
$$\sigma(x_{h_1},x_{j_1})\ge r+\delta>r.$$
Let $n_1$ be the smallest number with $n_1 > h_1 $ and 
$$\sigma( x_{h_1},x_{n_1}) \ge r+\delta.$$
Note 
$$\sigma(x_{h_1},x_{(n_1-1)})<r+\delta.$$
There exist $j_2>h_2>n_1$ such that $$\sigma(x_{h_2},x_{j_2})\ge r+\delta>r.$$
Let $n_2$ be the smallest number with $n_2 > m_2 $ and  
$$\sigma( x_{h_2},x_{n_2}) \ge r+\delta.$$
Then
$$\sigma(x_{h_2},x_{(n_2-1)})<r+\delta.$$
\indent Continuing this process, we build two increasing sequences in $\mathbb{N}$, $\{h_k\}_{k\in\mathbb{N}}$ and $\{n_k\}_{k\in\mathbb{N}}$ such that for all  $k$,
$$\sigma(x_{h_k},x_{(n_k-1)})<r+\delta\le \sigma(x_{h_k},x_{n_k})\hspace{6ex}(\triangledown)$$
and
$$m_{x_{h_k},x_{n_k}}=\sigma(x_{n_k},x_{n_k}).$$
For all  $k$, denote $s_k=\sigma(x_{h_k},x_{n_k})$. We should note that $$\sigma(x_{(h_k+1)},x_{(n_k+1)})\ne s_{k+1}=\sigma(x_{h_{(k+1)}},x_{n_{(k+1)}})$$
but rather, from \Cref{D7b},
$$\sigma(x_{(h_k+1)},x_{(n_k+1)})\le \sigma(x_{h_k},x_{n_k})-\varphi (\sigma(x_{h_k},x_{n_k}))=s_k-\varphi(s_k).\hspace{6ex}(\otimes)$$
Therefore, for all $k>N$ ($N$ defined in the beginning of step 2),

 $$r+\delta \le s_k=\sigma(x_{h_k},x_{n_k})=[\sigma(x_{h_k},x_{n_k})-m_{x_{h_k},x_{n_k}}]+m_{x_{h_k},x_{n_k}}$$
by ($\sigma-$ inq)
$$\le [\sigma(x_{h_k},x_{(n_k-1)})-m_{x_{h_k},x_{n_k-1}}+\sigma(x_{(n_k-1)},x_{n_k})-m_{x_{(n_k-1)},x_{n_k}}]+m_{x_{h_k},x_{n_k}}$$
by $(\triangledown)$ and step 1
$$\le r+\delta -r+t_{(n_k-1)}-r+m_{x_{h_k},x_{n_k}}.$$
Hence, taking $k \to +\infty$ by $(\bar\bigcirc)$ and  step 1 
$$r+\delta \le \lim_{k \to +\infty}s_k\le r+\delta-r+r-r+r=r+\delta.$$
We just showed that there exists $\delta >0$ such that
$$\lim_{k \to +\infty}s_k=r+\delta.$$
The next step is to prove that $\delta=0$ giving us our contradiction.
By applying ($\sigma-$inq) we get
$$s_k=\sigma(x_{h_k},x_{n_k})=[\sigma(x_{h_k},x_{n_k})-m_{x_{h_k},x_{n_k}}]+m_{x_{h_k},x_{n_k}}$$
$$\le [\sigma(x_{h_k},x_{(n_k+1)})-m_{x_{h_k},x_{(n_k+1)}}+\sigma(x_{(n_k+1)},x_{n_k})-m_{x_{(n_k+1)},x_{n_k}}]+m_{x_{h_k},x_{n_k}}$$
$$\le \sigma(x_{h_k},x_{(n_k+1)})-m_{x_{h_k},x_{(n_k+1)}}+t_{n_k}-r+m_{x_{h_k},x_{n_k}}.$$
Using ($\sigma-$inq) again on $\sigma(x_{h_k},x_{(n_k+1)})-m_{x_{h_k},x_{(n_k+1)}}$ we get
$$s_k \le \sigma(x_{h_k},x_{(h_k+1)})-m_{x_{h_k},x_{(h_k+1)}}+\sigma(x_{(h_k+1)},x_{(n_k+1)})-m_{x_{(h_k+1)},x_{(n_k+1)}}$$
$$+t_{n_k}-r+m_{x_{h_k},x_{n_k}}.$$
Therefore, by $(\otimes)$ and step 1,
$$s_k \le t_{h_k}-r+s_{k}-\varphi(s_k)-r+t_{n_k}-r+m_{x_{h_k},x_{n_k}}$$
i.e.
$$0 \le \varphi(s_k) \le t_{h_k}-r+-r+t_{n_k}-r+m_{x_{h_k},x_{n_k}}.$$
Taking the limit as $k \to +\infty$ we get
$$0 \le \lim_{k\to \infty} \varphi(s_k) \le r-r-r+r-r+r=0.$$
Hence,  and since $\varphi$ is non-decreasing with $r+\delta\le s_k$,
$$0\le \varphi(r+\delta)\le\lim_{k \to +\infty}\varphi(s_k)=0$$
i.e. $r+\delta =r$ and, therefore, $\delta=0$, a clear contradiction. Therefore, the assumption considered at the beginning of step 2 is incorrect proving that $\{x_i\}_{i\in \mathbb{N}}$ is an $r-$Cauchy sequence.$\hspace{6ex}\square$
\section{Continuity and Non-Expansiveness.}\label{S7}
\begin{defn}\label{D6a} Let $(X,\sigma)$ be an $M-$metric space with $x_o \in X$. We say that a function $f:X\to X$ is \textbf{\uline{weakly orbitally continuous at $x_o$}} iff 
\begin{center} if $a$ is the special limit of $\{f^i(x_o)\}_{i\in \mathbb{N}}$ then $f(a)$ is a limit of $\{f^i(x_o)\}_{i\in \mathbb{N}}$.
\end{center}
\end{defn}
\begin{rmk}\label{R6a} Notice that $f(a)$ is not required to be a special limit of $\{f^i(x_o)\}_{i\in \mathbb{N}}$, but rather only required to be its limit.
\end{rmk}
\begin{lem}\label{L6a} Let $(X,\sigma)$ be an $M-$metric space with $x_o \in X$. Let $f:X \to X$ be a weakly orbitally continuous function at $x_o$. If $a$ is a special limit of $\{f^i(x_o)\}_{i\in \mathbb{N}}$ then
$$m_{a,f(a)} =\sigma(a,a)\le \sigma(a,f(a))\le \sigma(f(a),f(a)).$$
\end{lem}
\textbf{Proof:} Denote for all natural numbers $i$, $x_i=f^i(x_o)$. Since $a$ is a special limit of  $\{f^i(x_o)\}_{i\in \mathbb{N}}$ and $f$ is weakly orbitally continuous at $x_o$, then $f(a)$ is a limit of $\{f^i(x_o)\}_{i\in \mathbb{N}}$. Therefore, by \Cref{L5b},
$$\sigma(a,a)\le\sigma(f(a),f(a)) $$
and, hence, by ($\sigma-$lbnd) 
$$m_{a,f(a)} =\sigma(a,a)\le \sigma(a,f(a)).$$
Furthermore, for all $i$,
$$\sigma(a,f(a))=[ \sigma(a,f(a))-m_{a,f(a)}]+\sigma(a,a)$$
by ($\sigma-$inq)
$$\le [\sigma(a,x_{i})-m_{a,x_i}+\sigma(x_i,f(a))-m_{x_i,f(a)}]+\sigma(a,a)$$
by adding and subtracting $\sigma(x_i,x_{i})$,
$$=\sigma(a,x_{i})-m_{a,x_i}+\uline{\uline{\sigma(x_i,f(a))+\sigma(x_i,x_{i})-m_{x_i,f(a)}}}-\sigma(x_i,x_{i})+\sigma(a,a)$$
taking the limit as $i \to +\infty$
$$=\sigma(a,a)-\sigma(a,a)+\uline{\uline{\sigma(f(a),f(a))}}-\sigma(a,a)+\sigma(a,a)=\sigma(f(a),f(a)).\hspace{6ex}\square$$
\indent In \Cref{L6a}, and by ($\sigma-$sep),  for the special limit $a$ to be a fixed point of $f$, we need $\sigma(f(a),f(a))\le\sigma(a,a)$. This can be obtained in various ways. In this paper we discuss two: The first is non-expansiveness, the second is the space having $\sigma(f(a),f(a))$ as a lower bound.
 \begin{defn}\label{D6b} Let $(X,\sigma)$ be an $M-$metric space. Let $f:X\to X$ be a function on $X$. We say that $f$ is \textbf{\uline{non-expansive}}  if and only if for all $x,y \in X$,
$$\sigma(f(x),f(y))\le\sigma(x,y).$$
\end{defn}
\begin{defn}\label{D6c} Let $(X,\sigma)$ be an $M-$metric space and $r_o$ a real number. We say on $(X,\sigma)$ is \uline{\textbf{bounded below by $r_o$}} if and only if for all $x,y \in X$,
$$r_o\le\sigma(x,y).$$
\end{defn}
\section{Main Theorem and Corollaries}\label{S8}
 We now present our main theorem. The reader will notice that we tried, as much as possible, to state it in its most general form.
 \begin{thm}
\label{T8a}
Let $(X,\sigma)$ be an  $M-$metric space with $x_o \in X$. 
Let $f:X\to X$ be a function such that $f$ is $r-$Cauchy at $x_0$ with special limit $a\in X$. Further assume at least one of the following conditions holds:
\begin{enumerate}
\item $f$ is weakly orbitally continuous at $x_0$ and  non-expansive.
\item $f$ is weakly orbitally continuous at $x_0$ and $(X, \sigma)$ is bounded below by $\sigma(f(a), f(a))$.
\item $f$ is non-expansive and $(X, \sigma)$ is bounded below by $\sigma(a, a)$.
\end{enumerate}
Then, $a$ is a fixed point of $f$.
\end{thm}
\textbf{Proof:}  
In (1) and (2), since $f$ is weakly orbitally continuous at $x_o$ then by \Cref{L6a}
$$\sigma(a,a)\le \sigma(a,f(a))\le \sigma(f(a),f(a)).$$
Both (1) $f$ is non-expansive and (2) $(X,\sigma)$ is bounded below by  $\sigma(f(a),f(a))$, assert $\sigma(f(a),f(a))\le\sigma(a,a)$. Therefore,
$$\sigma(a,a)= \sigma(a,f(a))= \sigma(f(a),f(a))$$
and, hence, by ($\sigma-$sep) $f(a)=a.$
As for (3), $f$ is non-expansive and  $(X,\sigma)$ is bounded below by  $\sigma(a,a)$ assert that $\sigma(f(a),f(a))=\sigma(a,a)$ and, hence, for all $i$,
$$\lim_{i\to+\infty}m_{f(a),f(x_i)}=\lim_{i\to+\infty}m_{a,x_i}=\sigma(a,a)=\sigma(f(a),f(a))=m_{f(a),a}.\hspace{6ex}(\ominus)$$
Therefore by ($\sigma-$lbnd), 
$$\sigma(a,a)=m_{a,f(a)}\le\sigma (a,f(a)).$$
By ($\sigma-$inq)  for all $i$,
$$\sigma(f(a),a)- m_{f(a),a}\le\sigma(f(a),f(x_i))-m_{f(a),f(x_i)}+\sigma(a,f(x_i))-m_{a,f(x_i)}.$$
Hence, by non-expansiveness,
$$\sigma(f(a),a)\le  m_{f(a),a}+\sigma(a,x_i)-m_{f(a),f(x_i)}+\sigma(a,f(x_i))-m_{a,f(x_i).}$$
Using $(\ominus)$ and by taking $i \to +\infty$ we get
$$\sigma(f(a),a)\le\sigma(a,a)+\sigma(a,a)-\sigma(a,a)+\sigma(a,a)-\sigma(a,a)=\sigma(a,a).$$
Therefore, by ($\sigma-$sep) $a=f(a)\hspace{6ex}\square$ 
\\\indent  Both  \Cref{L7a} and \Cref{L7b} assert that under their respective conditions, $f$ is an $r-$Cauchy sequence. Unfortunately, \Cref{T8a} presents us with the conditions to obtain a fixed point without guaranteeing its uniqueness. 
We now  present a valid proof of [\cite{Asa2014}: Theorem 3.1] and [\cite{Asa2014}: Theorem 3.2] in \Cref{C8a}, \Cref{C8b} respectively. \\\indent  We remind our reader that the definition presented for the $M-$metric in \cite{Asa2014} restricts $\sigma$ to non-negative values. This section is presented with that premise in mind.
\begin{cor}\label{C8a} Let $(X,\sigma)$ be a complete  $M-$metric space. Let $f:X \to X$ be a continuous function satisfying the following condition: There exists $0\le k < 1$ such that for all $x, y \in X$
$$0\le \sigma(f(x),f(y))\le k \sigma(x,y).\hspace{6ex} (\circ)$$
Then, $f$ has a unique fixed point.
\end{cor}
\textbf{Proof:} Consider any arbitrary $x_o \in X$.  The function   $f$ is  $\varphi_0-$contractive at $x_o$ ($\varphi_r$ with $r=0$) where $\varphi(t)=(1-k)t$. Hence, using \Cref{L7b}, $f$ is a $0-$Cauchy function at $x_o$. Since $(X,\sigma)$ is complete, let $a$ be the special limit of  $\{f^i(x_o)\}_{i\in \mathbb{N}}$. Since $f$ is continuous, then $f$ is weakly orbitally continuous at $x_o$. Additionally,  $(\circ)$ also asserts that $f$ is non-expansive. Therefore, by \Cref{T8a}(1), the special limit $a$ is a fixed point. Now to prove uniqueness. Assume that $a$ and $b$ are both fixed points of $f$. Hence, by $(\circ)$, 
$$\sigma(a,a)=\sigma(f(a),f(a))\le k \sigma(a,a)<\sigma(a,a),$$
$$\sigma(b,b)=\sigma(f(b),f(b))\le k \sigma(b,b)<\sigma(b,b)$$
and
$$\sigma(a,b)=\sigma(f(a),f(b))\le k \sigma(a,b)<\sigma(a,b).$$
Therefore, 
$$ \sigma (a,a)=\sigma(b,b)=\sigma(a,b)=0$$
and, hence, by ($\sigma-$sep) $a=b$. $\hspace{6ex}\square$
\begin{cor}\label{C8b} Let $(X,\sigma)$ be a complete  $M-$metric space. Let $f:X \to X$ be a continuous function satisfying the following condition: There exists $0\le k < \frac{1}{2}$ such that for all $x, y \in X$
$$0\le \sigma(f(x),f(y))\le k [\sigma(x,f(x))+\sigma(y,f(y))].\hspace{6ex} (\triangle)$$
Then, $f$ has a unique fixed point.
\end{cor}
\textbf{Proof:}  Consider any arbitrary $x_o \in X$ and denote  $x_i =f^i(x_o)$. We first show that $f$ is a $c_0-$contraction at $x_o$ ($c_r$ with $r=0$) and $c=2k>\frac{k}{1-k}$. \\By ($\triangle$) we get for every $i$,
$$\sigma(x_{i+1},x_{i+1})\le 2k\sigma(x_{i+1},x_{i}).$$
Moreover, for every $i$,
$$\sigma(x_{i+2},x_{i+1})\le k(\sigma(x_{i+1},x_{i+2})+\sigma(x_{i},x_{i+1}))$$
i.e.
$$\sigma(x_{i+2},x_{i+1})\le \frac{k}{1-k}\sigma(x_{i+1},x_{i})<c\sigma(x_{i+1},x_{i}).$$
Hence,
$$0\le \sigma(x_{i+2},x_{i+2})<\sigma(x_{i+2},x_{i+1})\le c^{i+1}\sigma(x_{1},x_{o})$$
and by \Cref{L7a}, $f$ is a $0-$Cauchy function at $x_o$. Since $(X,\sigma)$ is complete, $\{x_i\}_{i\in\mathbb{N}}$ has a special limit $a$. Hence, by \Cref{D5c}, $\sigma(a,a)=0$.
\\\indent  The function $f$ is continuous and, hence, weekly orbitally continuous at $x_o$. Therefore, by  \Cref{L6a}, we have
$$m_{a,f(a)} =\sigma(a,a)\le \sigma(a,f(a))\le \sigma(f(a),f(a)).$$ 
Additionally, by $(\triangle)$,
$$\sigma(f(a),f(a))\le 2k \sigma(a,f(a))\le 2k\sigma(f(a),f(a)).$$
Hence, $\sigma(f(a),f(a))=0$ completing the requirement for \Cref{T8a}(2). As for uniqueness,
let $a$ and $b$ are both fixed points of $f$. Hence, by $(\triangle)$, 
$$\sigma(a,b)=\sigma(f(a),f(b))\le 2k\sigma(a,b)$$
$$\sigma(a,a)=\sigma(f(a),f(a))\le 2k\sigma(a,a)$$
and, similarly,
$$\sigma(b,b)=\sigma(f(b),f(b))\le 2k\sigma(b,b).$$
Therefore,
$$\sigma(a,a)=\sigma(b,b)=\sigma(a,b)=0$$
and, by ($\sigma-$sep), $a=b.\hspace{6ex}\square$

\bibliographystyle{amsplain}
\bibliography{Generalizedmetrics}


\end{document}